\documentclass{amsart}
\usepackage{amsfonts}
\usepackage{amsmath}
\usepackage{amssymb}
\usepackage{amsthm}
\usepackage{url}

\newtheorem{theorem}{Theorem}[section]

\newtheorem{corollary}[theorem]{Corollary}

\newtheorem{lemma}[theorem]{Lemma}

\theoremstyle{remark}
\newtheorem*{remark}{Remark}

\newcommand{\A}{\mathbb{A}}
\newcommand{\GF}{\mathbb{F}}

\newcommand{\calA}{\mathcal{A}}
\newcommand{\calB}{\mathcal{B}}
\newcommand{\calC}{\mathcal{C}}
\newcommand{\calM}{\mathcal{M}}
\newcommand{\calS}{\mathcal{S}}

\newcommand{\GFbar}{\overline{\mathbb{F}}}
\newcommand{\bbar}{\overline{b}}
\newcommand{\cbar}{\overline{c}}
\newcommand{\kbar}{\overline{k}}
\newcommand{\epsbar}{\overline{\eps}}
\newcommand{\sigmabar}{\overline{\sigma}}
\newcommand{\taubar}{\overline{\tau}}

\newcommand{\eps}{\varepsilon}

\DeclareMathOperator{\Aut}{Aut}
\DeclareMathOperator{\divisor}{div}
\DeclareMathOperator{\End}{End}
\DeclareMathOperator{\Gal}{Gal}
\DeclareMathOperator{\Jac}{Jac}
\DeclareMathOperator\ord{ord}
\DeclareMathOperator{\PSL}{PSL}
\DeclareMathOperator{\Res}{Res}
\DeclareMathOperator{\Tr}{Tr}

\begin{document}

\title[Supersingular genus-two curves]
{Supersingular genus-two curves \\ over fields of characteristic three}

\author{Everett W.~Howe}
\address{Center for Communications Research, 
         4320 Westerra Court, 
         San Diego, CA 92121-1967, USA.} 
\email{however@alumni.caltech.edu}
\urladdr{http://www.alumni.caltech.edu/\~{}however/}

\date{18 April 2006}



\begin{abstract}
Let $C$ be a supersingular genus-$2$ curve over an algebraically closed
field of characteristic~$3$.  We show that if $C$ is not isomorphic to
the curve $y^2 = x^5 + 1$ then up to isomorphism there are exactly $20$
degree-$3$ maps $\varphi$ from $C$ to the elliptic curve $E$ with 
$j$-invariant~$0$.  We study the coarse moduli space of triples 
$(C,E,\varphi)$, paying particular attention to questions of rationality.
The results we obtain allow us to determine, for every finite field $k$ of 
characteristic~$3$, the polynomials that occur as Weil polynomials of 
supersingular genus-$2$ curves over~$k$.
\end{abstract}

\maketitle

\section{Introduction}
\label{S:Introduction}

Over an algebraically closed field of characteristic~$3$, every
supersingular genus-$2$ curve $C$ has a degree-$3$ map to the 
supersingular elliptic curve $E$ with $j$-invariant~$0$.  We study the 
coarse moduli space $\calA$ of such maps $C\to E$, paying special
attention to questions of rationality.  We show that with a single
exception, every supersingular genus-$2$ curve over an 
algebraically-closed field of characteristic $3$ has $20$ non-isomorphic
maps to the elliptic curve with $j$-invariant~$0$, and we provide
explicit equations for these maps. 

Our motivation for studying these maps is the question of determining
which isogeny classes of abelian surfaces over a finite field contain 
Jacobians.  Much has been written on this question
\cite{
AdlemanHuang:LNM, 
Howe:TAMS, 
Howe:JAG, 
Howe:Texel, 
Howe:CM, 
HMNR,
Maisner:Thesis, 
MaisnerNart:EM, 
MaisnerNart:CJM, 
McGuireVoloch:PAMS, 
Ruck:CM, 
Serre:notes}, and together the works just cited answer the question
for all simple non-supersingular isogeny classes.  In recent joint
work~\cite{HNR}, Enric Nart, Christophe Ritzenthaler, and the author
developed techniques to address the non-simple isogeny classes and 
the supersingular isogeny classes, thus completing the determination
of the isogeny classes of abelian surfaces over finite fields that
contain Jacobians.  However, to avoid various special cases, the 
arguments in~\cite{HNR} concerning supersingular isogeny classes assume
that the characteristic of the base field is larger than~$3$.  The 
characteristic~$2$ case is dealt with in~\cite{MaisnerNart:CJM}.  
In this paper, we use our results on the moduli space $\calA$ to answer
the question for supersingular isogeny classes in characteristic~$3$.
We obtain the following theorem:
\begin{theorem}
\label{T:intro}
Let $q = 3^d$ be a power of $3$.  If $d$ is odd, then the polynomials
that occur as the Weil polynomials of supersingular genus-$2$ curves
over $\GF_q$ are{\rm:}
\begin{enumerate}
\item $(x^2 + q)(x^2 - sx + q)$ for all $s\in\{\pm\sqrt{3q}\}${\rm;}
\item $(x^2 + q)^2$, if $q > 3${\rm;}
\item $x^4 + q^2${\rm;}
\item $x^4 + qx^2 + q^2${\rm;}
\item $x^4 - 2qx^2 + q^2$, if $q > 3$.
\end{enumerate}
If $d$ is even, the polynomials that occur as the Weil polynomials of 
supersingular genus-$2$ curves over $\GF_q$ are{\rm:}
\begin{enumerate}
\item $(x^2 - 2sx + q)(x^2 + sx + q)$ for all $s\in\{\pm\sqrt{q}\}${\rm;}
\item $(x^2 - sx +  q)^2$ for all $s\in\{0,\pm\sqrt{q}\}${\rm;}
\item $(x^2 - 2sx +  q)^2$ for all $s\in\{\pm\sqrt{q}\}$, if $q > 9${\rm;}
\item $x^4 + q^2${\rm;}
\item $x^4 - sx^3 + qx^2 - sq x+ q^2$ for all $s\in\{\pm\sqrt{q}\}$.
\end{enumerate}
\end{theorem}

To prove these results, we show that the moduli space $\calA$ is 
isomorphic to the coarse moduli space $\calB$ of pairs $(C,G)$, where 
$C$ is a supersingular genus-$2$ curve in characteristic $3$ and $G$ is 
an order-$4$ subgroup of $(\Jac C)[2]$ that is not isotropic with respect
to the Weil pairing.  We give explicit models for these moduli spaces, as
well as for the maps from these spaces to the moduli space of supersingular
genus-$2$ curves in characteristic~$3$.

In Section~\ref{S:basic3} we provide some basic facts about supersingular
genus-$2$ curves in characteristic~$3$; in particular, we show that the
coarse moduli space $\calS$ of such curves is isomorphic to the affine
line.  In Section~\ref{S:covers} we classify the degree-$3$ maps from 
supersingular genus-$2$ curves to the $j$-invariant $0$ elliptic curve 
in characteristic~$3$.  We show that the coarse moduli space $\calA$ of 
such maps is isomorphic to the affine line with one point removed, and 
that the natural map $\calA\to\calS$ has degree~$20$.  
In Section~\ref{S:SSEC} we give complete lists of the isomorphism 
classes of supersingular elliptic curves over every finite field of 
characteristic~$3$, and we calculate the Weil polynomial of each
isomorphism class.  Finally, in Section~\ref{S:WeilPolys3} we use the 
results of the earlier sections to determine which isogeny classes of
supersingular abelian surfaces in characteristic~$3$ do not contain
Jacobians, thereby proving Theorem~\ref{T:intro}.

\subsubsection*{Notation}
If $X$ is a scheme over a field $K$, and if $L$ is an extension field
of~$K$, we let $X_L$ denote the base extension of $X$ to $L$.

\subsubsection*{Acknowledgments}
The author is grateful to Enric Nart and Christophe Ritzenthaler
for encouragement and helpful discussions.

\section{Basic facts about supersingular genus-two curves in characteristic
three.}
\label{S:basic3}

In this section we determine the locus of supersingular genus-$2$ curves
in characteristic~$3$, viewed as a subvariety of the coarse moduli space
of genus-$2$ curves.

If $k$ is a field, we set $S_k = \A_k^5 \setminus\{[0,0,0,0,0]\}$
and we let $W$ be the orbit space of $S_k$ under the action of $k^*$
defined by 
$$\lambda[x_1,x_2,x_3,x_4,x_5] = 
 [\lambda^2 x_1, \lambda^4 x_2, \lambda^6 x_3, \lambda^8 x_4, \lambda^{10} x_5],$$
so that $W$ is a weighted projective space.  We denote the orbit of
$[x_1,x_2,x_3,x_4,x_5]$ by $[x_1:x_2:x_3:x_4:x_5]$.  Igusa~\cite{Igusa:AM} 
associated to every genus-$2$ curve over $k$ an element 
$[J_2:J_4:J_6:J_8:J_{10}]$ of $W$ such that $J_{10}\neq 0$ and 
$4J_8 = J_2 J_6 - J_4^2$.  This element is the 
\emph{vector of Igusa invariants} of the curve.  Igusa showed that over
an algebraically closed field, every element of $W$ with $J_{10}\neq 0$
and $4J_8 = J_2 J_6 - J_4^2$ comes from a genus-$2$ curve, and  two 
genus-$2$ curves have the same vector of Igusa invariants if and only
if they are isomorphic to one another.  (Note, however, that a point of
this variety that is rational over a subfield of $k$ does not always come 
from a curve defined over that field.)  Thus, Igusa determined the coarse 
moduli space $\calM_2$ of genus-$2$ curves.

\begin{theorem}
\label{T:SSfacts1}
Let $k$ be a field of characteristic $3$.
\begin{itemize}
\item[(a)]
Let $C$ be a genus-$2$ curve over $k$ and let $[J_2:J_4:J_6:J_8:J_{10}]$ be
its vector of Igusa invariants.  Then $C$ is supersingular if and only 
if $J_2 = J_4 = J_8 = 0$.
\item[(b)]
Let $A$ and $B$ be elements of $k$, with $B\neq 0$.  Then the curve
$$y^2 = x^6 + Ax^3 + Bx + A^2$$
is a supersingular curve of genus~$2$, and its vector of Igusa invariants
is equal to $[0:0:A:0:B]$.
\end{itemize}
\end{theorem}

\begin{proof}
First we dispose of an annoying special case.  Suppose $k=\GF_3$ and $C$
is the curve $y^2 = x^5 - x$.  Then $C$ is not supersingular (because 
its Weil polynomial is $t^4 - 2t^2 + 9$) and its vector of Igusa invariants
is $[1:0:-1:-1:-1]$.

If $k$ is larger than $\GF_3$, or if $C$ is any curve other than
$y^2 = x^5 - x$, then $C$ has a model over $k$ of the form
$d y^2 = x^6 + c_5 x^5 + c_4 x^4 + c_3 x^3 + c_2 x^2 + c_1 x + c_0,$
where the polynomial in $x$ is separable.  Now we use the criterion\
for supersingularity given in~\cite[Lemma~E, p.~387]{Yui:JA}:
Let $M$ and $M^{(3)}$ be the matrices
$$ M = \left[\begin{matrix}c_2&c_1\\c_5&c_4\end{matrix}\right] \qquad
 M^{(3)} = \left[\begin{matrix}c_2^3&c_1^3\\c_5^3&c_4^3\end{matrix}\right].$$
Then $C$ is supersingular if and only if $M M^{(3)} = 0 $, and the Jacobian
of $C$ is isomorphic to a product of supersingular elliptic curves
(that is, the Jacobian is \emph{superspecial}) if and only if $M = 0$. 
(Actually, to apply~\cite[Lemma~E, p.~387]{Yui:JA} we should write our
curve as $y^2 = \text{(quintic)}$.  However, the proof in~\cite{Yui:JA}
works for curves of the form $y^2 = \text{(sextic)}$ as well.)

We note in passing that $M\neq 0$, because a polynomial of the form 
$x^6 + c_3 x^3 + c_0$ is not separable.  This shows that there are no 
superspecial genus-$2$ curves in characteristic~$3$, as we know
already from~\cite[Prop.~3.1]{IKO:CM}.

Note that 
$$M M^{(3)} = \left[\begin{matrix}
c_2^4 + c_1 c_5^3     & c_2 c_1^3 + c_1 c_4^3 \\
c_5 c_2^3 + c_4 c_5^3 & c_5 c_1^3 + c_4^4
\end{matrix}\right].$$
Suppose the curve is supersingular, so that this matrix is the zero
matrix.  If $c_5 = 0$ then we must have $c_4 = c_2 = 0$, so $C$ can be
written in the form $d y^2 = x^6 + c_3 x^3 + c_1 x + c_0.$
On the other hand, if $c_5$ is nonzero then we can translate $x$ by
a constant in order to get a new equation for $C$ of the form
$d y^2 = x^6 + c_5 x^5 + c_3 x^3 + c_2 x^2 + c_1 x + c_0$.
Then we find that $c_2 = c_1 = 0$, so $C$ is of the form
$d y^2 = x^6 + c_5 x^5 + c_3 x^3 +  c_0$.
If we replace $x$ with $1/x$ and rescale $y$, we find that 
we can again write $C$ in the form
$d y^2 = x^6 + c_3 x^3 + c_1 x + c_0.$

One can calculate that the Igusa invariants of the curve 
$d y^2 = x^6 + c_3 x^3 + c_1 x + c_0$ are
$[0: 0: c_3^6 - c_3 c_1^3 - c_0^3: 0: -c_1^6]$.
(The computer algebra system Magma~\cite{Magma} provides facilities for
calculating Igusa invariants.)  This proves the `only if' part of 
statement~(a).

If we take $c_3 = A$, $c_1 = B$, $c_0 = A^2$, and $d=1$, we find that
the Igusa invariants are  $[0: 0:  - A B^3: 0: -B^6] = [0:0:A:0:B]$.
Thus, any curve with $J_2 = J_4 = J_8 = 0$ is geometrically isomorphic
to one of the form $y^2 = x^6 + Ax^3 + Bx + A^2$, and the criterion 
from~\cite{Yui:JA} shows that this curve is supersingular.  This proves
statement (b) and the `if' part of statement~(a).
\end{proof}

Given a supersingular genus-$2$ curve $C$ over a field of
characteristic~$3$, we let $I_C = J_6^5/J_{10}^3$, where
$[0:0:J_6:0:J_{10}]$ is the vector of Igusa invariants of $C$.

\begin{theorem}
\label{T:SSfacts2}
Let $k$ be a field of characteristic~$3$.
\begin{itemize}
\item[(a)]
The map
$$\left\{\vcenter{\hsize=5.5cm \noindent
\hfil$k$-isomorphism classes of\hfil\break
supersingular genus-$2$ curves over $k$}\right\}
\to k$$
that sends $C$ to $I_C$ is surjective, and is bijective if $k$ 
is algebraically closed.
\item[(b)]
If $c$ is a nonzero element of $k$, then the curve $C$ defined by
$$y^2 = x^6 + c^2 x^3 + c^3 x + c^4$$
is supersingular and satisfies $I_C = c$.  The curve $C$ defined by
$y^2 = x^5 + 1$ is supersingular and satisfies $I_C = 0$.
\item[(c)]
Let $C$ be a supersingular genus-$2$ curve over $k$.  The geometric 
automorphism group of $C$ has order $2$ if $I_C\neq 0$, and order $10$
if $I_C = 0$.
\end{itemize}
\end{theorem}

\begin{proof}
Given any nonzero $c$, let $C$ be the curve 
$y^2 = x^6 + c^2 x^3 + c^3 x + c^4$.  It is easy to check that
$I_C = c$.  Likewise, it is easy to check that the Igusa invariants 
of the curve $C$ defined by $y^2 = x^5 + 1$  are $[0:0:0:0:1]$, so we
have $I_C = 0$.  This proves statement (b) and the surjectivity claim 
of statement~(a).

Suppose $C$ and $C'$ are supersingular genus-$2$ curves, with 
Igusa invariants $[0:0:A:0:B]$ and $[0:0:A':0:B']$, respectively. 
Then $I_C = A^5/B^3$ and $I_{C'}=(A')^5 / (B')^3$.  Let $\lambda$ be 
a $10$th root of $B'/B$.  If $I_C = I_{C'}$ then we have 
$(A')^5 = \lambda^{30} A^5$, so $A' = \zeta \lambda^6 A$ for some
$5$-th root of unity $\zeta$.  If we replace $\lambda$ with 
$\lambda/\zeta$, then we have $B' = \lambda^{10} B$ and 
$A' = \lambda^6 A$, so $[0:0:A:0:B] = [0:0:A':0:B']$, and $C$ and 
$C'$ are geometrically isomorphic to one another.  This proves the
bijectivity claim of statement (a).

Finally, we note that Igusa~\cite[\S8]{Igusa:AM} calculated the 
automorphism groups of all genus-$2$ curves in every characteristic. 
In characteristic~$3$, the only curve with automorphism group larger 
than $\{\pm 1\}$ is the curve with Igusa invariants $[0:0:0:0:1]$, 
whose automorphism group is cyclic of order $10$.  This proves 
statement~(c).
\end{proof}

\begin{remark}
Results similar to Theorem~\ref{T:SSfacts1} can be found 
in~\cite[\S6]{Zarhin:BSMF} (see also the corrections
in~\cite[\S5]{Zarhin:AGCT9}).
\end{remark}

\begin{remark}
Let $\calS$ denote the coarse moduli space of supersingular genus-$2$
curves in characteristic~$3$.  We have shown that the invariant $I_C$ 
gives an isomorphism from $\calS$ to the affine line $\A^1$.
\end{remark}

\section{Supersingular genus-two curves as triple covers of elliptic curves.}
\label{S:covers}

In this section we classify all of the degree-$3$ maps 
$\varphi\colon C\to E$ over a given base field of characteristic~$3$,
where $C$ is a genus-$2$ curve and $E$ is a supersingular elliptic 
curve. We start by writing down a family of such maps.

Let $k$ be a field of characteristic $3$, and let $b$ and $c$ be nonzero
elements of $k$.  Let $E_{b,c}$ be the elliptic curve
$$y^2 = x^3 - bx + c$$
and let $C_{b,c}$ be the genus-$2$ curve
$$w^2 = c (v^3 - bv^2 - b^2v + b^3 - c^2) (v^3 + bv^2 - b^2v - b^3 - c^2).$$

\begin{lemma}
\label{L:SSsplitting}
There are degree-$3$ maps 
$$\varphi_{b,c}\colon C_{b,c}\to E_{b,c}
       \text{\qquad and\qquad}
  \varphi'_{b,c}\colon C_{b,c}\to E_{-b,c}$$
that induce an isogeny $\Jac C_{b,c} \to E_{b,c}\times E_{-b,c}$.
\end{lemma}

\begin{proof}
Let 
\begin{align*}
x &= -bc (v-b) / (v^3 - bv^2 - b^2 v + b^3 - c^2)\\
z &= -w /(v^3 - bv^2 - b^2 v + b^3 - c^2)\\
y &= (z^3 + bxz) / c.
\end{align*}
One can check that then $y^2 = x^3 - bx + c$, so these equations
define a degree-$3$ map $\varphi_{b,c}\colon C_{b,c}\to E_{b,c}$.
Replacing $b$ with $-b$ in the definition of $\varphi_{b,c}$
gives a degree-$3$ map $\varphi'_{b,c}\colon C_{b,c}\to E_{-b,c}$.
It is an easy exercise to pull back the invariant differentials
of $E_{b,c}$ and $E_{-b,c}$ to $C_{b,c}$ and to show that the
pullbacks are linearly independent; it follows that the induced
map $\Jac C_{b,c} \to E_{b,c}\times E_{-b,c}$ is an isogeny.
\end{proof}

\begin{theorem}
\label{T:SSnormalform}
Suppose $C$ is a curve over a field $k$ of characteristic $3$, and
suppose $\varphi\colon C\to E$ is a degree-$3$ map from $C$ to a
supersingular elliptic curve $E$ over~$k$.  Then there exist nonzero
$b,c\in k$ such that the cover $\varphi\colon C\to E$ is isomorphic
$\varphi_{b,c}\colon C_{b,c}\to E_{b,c}$.
\end{theorem}

\begin{proof}
The Riemann-Hurwitz formula shows that the degree of the different
of $\varphi$ is $2$, so either there are two geometric points of $C$
each contributing $1$ to the different, or there is one geometric point
contributing $2$ to the different.  A ramification analysis as 
in~\cite[\S7.2]{HoweLauter:AIF} shows that no point can contribute $2$
to the different; therefore, the triple cover $\varphi$ must have exactly
two ramification points $P$ and $Q$, both tame.\footnote{
  The analysis in~\cite[\S7.2]{HoweLauter:AIF} assumes that the base
  field is finite and that the cover has been written in a certain
  standard form.  However, for the ramification analysis all we
  need to know is that the cover can be put in the standard form
  locally everywhere, and this is true.
  
  We should note that there is an error in~\cite[\S7.2]{HoweLauter:AIF},
  starting at the second full paragraph on page 1717.  The error is in
  the statement that either $2\ord_P g \ge 3\ord_P f$ or 
  $\ord_P g \not\equiv 0\bmod 3$, except when $P = \infty$ and 
  $\ord_P g = -3$.  While this is true when $\ord_P f = 0$, it can 
  fail when $\ord_P f = 1$.  What \emph{is} true is that for every 
  $P\neq\infty$ for which $\ord_P f > 0$, there is a constant 
  $c_P\in\kbar$ such that either 
  $2\ord_P (g+c_P^3 - c_P f) \ge 3\ord_P f$ or 
  $\ord_P (g+c_P^3 - c_P f) \not\equiv 0\bmod 3$.  This change affects
  the statement in the final paragraph of~\cite[\S7.2]{HoweLauter:AIF}.
  See~\cite{HoweLauter:Corr} for more details.}
Let $k(C)$ be the function field of $C$ and let $k(E)$ be the function
field of $E$.  We can choose a generator $z$ for $k(C)$ over $k(E)$ whose
trace to $k(E)$ is $0$; then $z$ satisfies an equation $z^3 - fz = g$ 
for two functions $f$ and $g$ on $E$.  From~\cite[\S7.2]{HoweLauter:AIF}
we see that $f$ has odd order at a point of $E$ if and only if that
point contributes an odd amount to the different of $\varphi$; therefore
we must have
$$\divisor f = P + Q + 2D - 2\infty$$
where $D$ is a degree-$0$ divisor that is rational over~$k$. We can write
$$D =  - R + \infty + \divisor r$$
for some function $r$ and some point $R$ on $E$, both rational over~$k$.
Replacing $z$ with $z/r$, the function $f$ with $f/r^2$, and $g$ 
with $g/r^3$, we still have $z^3 - fz = g$, but now the divisor of 
$f$ is $P + Q - 2R$. (Note that $P$ and $Q$ must be distinct because 
there are two ramification points, and $R$ must be distinct from both
$P$ and $Q$, because otherwise the function $f$ would have a single zero
and a single pole.)  Now we can change co\"ordinates on $E$ so that $R$
is the point at infinity.  After this change of co\"ordinates, we have 
$P + Q = 0$ in the group law on $E$.  Let the equation for $E$ be 
$y^2 = x^3 - bx + c$, where $b\neq 0$.  Then the function $f$ must be of
the form $f = d(x-a)$ for some $a\in k$ and $d\in k^*$.  Shifting $x$ by 
the constant $a$ does not change the form of the equation for $E$, but now 
we may assume that $f = dx$ for some nonzero $d\in k$.  Since $P$ and $Q$
are distinct and both are zeros of $f$, we see that $c\neq 0$.

We know from \cite[\S7]{HoweLauter:AIF} that by translating $z$ by 
functions on $E$, we can modify $g$ so that it only has poles at the
poles of $f$; thus we may assume that the only pole of $g$ is at~$\infty$.
Furthermore, using the statement from the first full paragraph on 
page~1717 of~\cite{HoweLauter:AIF}, we see that we can modify $g$ so 
that it has at most a triple pole at~$\infty$.  Therefore, we may assume
that $g$ is a linear combination of $1$, $x$, and $y$.  By shifting $z$
by an appropriate constant to eliminate the multiple of $x$ appearing in
this linear combination, we may assume that $g$ is a linear combination
of $1$ and $y$.

Now we use the statement from the final paragraph 
of~\cite[\S7.2]{HoweLauter:AIF}, as corrected in the footnote above. 
We know that there are constants $c_P$ and $c_Q$ in $\kbar$ such that 
\begin{align*}
2\ord_P (g + c_P^3 - c_P f) &\ge 3\ord_P f = 3\\
\intertext{and}
2\ord_Q (g + c_Q^3 - c_Q f) &\ge 3\ord_Q f = 3.
\end{align*}
The functions on the left-hand sides of these inequalities have degree
at most $3$, so their orders at $P$ (and $Q$) are either $2$ or $3$. 
This means that the two lines
\begin{equation}
\label{EQ:lines1}
g + c_P^3 - c_P f = 0 \text{\qquad and\qquad}
g + c_Q^3 - c_Q f = 0
\end{equation}
must be the tangent lines to $E$ at $P$ and $Q$.

Let us write $P = (0, y_P)$ and $Q = (0, y_Q)$, with $y_Q = -y_P$.
Then the tangent lines to $E$ at $P$ and $Q$ are given by
\begin{equation}
\label{EQ:lines2}
y_P y - b x - y_P^2 = 0 \text{\qquad and\qquad}
y_Q y - b x - y_Q^2 = 0
\end{equation}
Note that the slopes of the lines defined by the these two equations are
negatives of one another, and the same is true of the $y$-intercepts of 
the two lines.  Applying the first of these two observations to the lines
defined in \eqref{EQ:lines1} shows that $c_Q = -c_P$.  Applying the second
observation then tells us that there can be no constant term appearing in
the function $g$; that is, $g = ey$ for some constant $e \in k$.  Finally,
by using this formula for $g$ and the fact that the two pairs of 
lines~\eqref{EQ:lines1} and~\eqref{EQ:lines2} are equal, we find that we
must have $b^3 e^2 + c^2 d^3 = 0$.  Since $b$, $c$, and $d$ are nonzero,
so is~$e$.

We have shown that the triple cover $C$ of $E$ can be written in the form 
\begin{equation}
\label{EQ:cover}
\begin{split}
y^2 &= x^3 - bx + c\\
z^3 - dxz &= ey
\end{split}
\end{equation}
where $b$, $c$, $d$, and $e$ are nonzero elements of $k$ that satisfy
$b^3 e^2 + c^2 d^3 = 0$.  It is convenient to rescale the variables in
these equations to get a more standardized form.  Let $r$ be a nonzero 
element of $k$, and replace $z$ with $rz$, $d$ with $d/r^2$, and $e$ 
with $e/r^3$.  Then the equations~\eqref{EQ:cover} still hold, but the 
ratio $c/e$ has been replaced with the ratio $cr^3/e$.  If we choose
$r = -(be)/(cd)$, then 
$$\frac{cr^3}{e} = \frac{-c b^3 e^3}{e c^3 d^3} 
                 = \frac{-b^3 e^2}{c^2 d^3} = 1.$$
In other words, we may scale our variables so that $e = c$, from which
it follows that $d = -b$.

We can find new equations defining $C$ by setting
\begin{align*}
v & = -b + (z^2 - c) /x\\
w & = - z ( v^3  - b v^2 - b^2 v + b^3 - c^2).
\end{align*}
A straightforward computation then shows that $v$ and $w$ satisfy
$$w^2 = c (v^3 - bv^2 - b^2v + b^3 - c^2)
          (v^3 + bv^2 - b^2v - b^3 - c^2),$$
and that we have
\begin{align*}
x &= -bc (v-b) / (v^3 - bv^2 - b^2 v + b^3 - c^2)\\
z &= -w /(v^3 - bv^2 - b^2 v + b^3 - c^2)\\
y &= (z^3 + bxz) / c.
\end{align*}
This shows that any given genus-$2$ triple cover of a supersingular
elliptic curve over a field $k$ of characteristic $3$ is isomorphic to 
a cover $\varphi_{b,c}\colon C_{b,c}\to E_{b,c}$ for some nonzero $b$ 
and $c$ in $k$.
\end{proof}


\begin{corollary}
\label{C:SScovers}
If $C$ is a genus-$2$ triple cover of a supersingular elliptic curve 
in characteristic~$3$, then $C$ is supersingular. \qed
\end{corollary}

Note that if $C$ is a triple cover of a supersingular elliptic curve, 
then the right-hand side of the equation for $C$ factors into a product
of two cubics.  The converse is true as well. 

\begin{theorem}
\label{T:cubicfactors}
Let $C$ be a supersingular genus-$2$ curve over a field $k$ of 
characteristic~$3$, and write $C$ in the form $y^2 = f$ for a sextic 
polynomial $f$.  If $f$ can be written as the product of two cubic 
factors, then $C$ is a triple cover of a supersingular elliptic curve.
\end{theorem}

\begin{proof}
Neither the hypothesis nor the conclusion of the theorem depends on
which model for $C$ we choose, so we may replace the given $f$ with
any sextic that defines the same curve.

First we note that some quadratic twist of $C$ can be written in the
form $y^2 = f$, where 
$$f = x^6 + A x^3 + B x + A^2$$
for some elements $A,B\in k$ with $B\neq 0$.  For curves with $I_C\neq 0$
this follows from Theorem~\ref{T:SSfacts2}(c); for curves with $I_C = 0$
this follows from explicitly writing down the twists of the curve
$y^2 = x^5 + 1$.

If we replace $x$ in this equation with $-B x$, and if we replace $y$ 
with $B^3 y$, we see that the new variables satisfy an equation of the
same form as above, only now the new $B$ has the property that $-B$ is 
a square.

Let $g_1$ and $g_2$ be monic cubics such that $f = g_1 g_2$.  Elementary
manipulations show that the coefficients of $x^2$ in $g_1$ and $g_2$ are
negatives of one another, and they are both nonzero.  Likewise, it is
easy to show that the sum of the coefficients of $x$ in $g_1$ and $g_2$
must be the square of the coefficient of $x^2$.  Thus we may write
\begin{align*}
g_1 &= x^3 + u x^2 + u^2 v x + w\\
g_2 &= x^3 - u x^2 + u^2 (1-v) x + t
\end{align*}
for some elements $u,v,w,t$ of $k$ with $u\neq 0$.  Setting $g_1 g_2$
equal to $f$, we find that we must have
\begin{align*}
t &= w + u^3(v^2-v)\\
\intertext{and}
w(v+1) &= -u^3 (v^4 - v^2 +1).
\end{align*}
The latter equality shows that $v\neq -1$.

For $i=1,2$ let $h_i$ be $g_i$ evaluated at $x + u(v+1)$.  If we let $c$
be a square root of $-B / u^2$ (which is possible because of the way we
adjusted $B$ at the beginning of the proof), then we have
\begin{align*}
h_1 &= x^3 + u x^2 - u^2 x + u^3 - c^2\\
h_2 &= x^3 - u x^2 - u^2 x - u^3 - c^2.
\end{align*}
Thus some quadratic twist of our curve $C$ is isomorphic to the curve
$C_{u,c}$.  But it is easy to see that a quadratic twist of a curve 
$C_{b,c}$ is isomorphic to a curve of the same form, so in any case, 
$C$ is a triple cover of an elliptic curve.
\end{proof}

Next we note that different $(b,c)$ pairs can give rise to isomorphic
covers $\varphi_{b,c}$.

\begin{lemma}
\label{L:IsomorphicCovers}
Let $b$, $c$, and $r$ be nonzero elements of $k$.  Then the cover
$$\varphi_{br^4,cr^6}\colon C_{br^4,cr^6}\to E_{br^4,cr^6}$$
is isomorphic to the cover $\varphi_{b,c}\colon C_{b,c}\to E_{b,c}$.
\end{lemma}

\begin{proof}
We can view the cover $\varphi_{b,c}\colon C_{b,c}\to E_{b,c}$ as being
given by the equations~\eqref{EQ:cover}, with $e = c$ and $d = -b$.
If we let $X = r^2 x$ and $Y = r^3 y$ and $Z = r^3 z$, then these
equations become
\begin{equation*}
\begin{split}
Y^2 &= X^3 - br^4X + cr^6\\
Z^3 + br^4 XZ &= cr^6 Y,
\end{split}
\end{equation*}
which are the equations defining the cover 
$$\varphi_{br^4,cr^6}\colon C_{br^4,cr^6}\to E_{br^4,cr^6}.$$
\end{proof}

Let $\calA$ denote the coarse moduli space of triples $(C,E,\varphi)$,
where $\varphi\colon C\to E$ is a degree-$3$ map from a supersingular
genus-$2$ curve in characteristic $3$ to a supersingular elliptic curve.
Let $\calB$ denote the coarse moduli space of pairs $(C,G)$, where $C$ 
is a supersingular genus-$2$ curve in characteristic $3$ and $G$ is an 
order-$4$ subgroup of $(\Jac C)[2]$ that is not isotropic with respect
to the Weil pairing.  There is a map $\kappa\colon\calA\to\calB$ that 
sends $(C,E,\varphi)$ to $(C,(\ker\varphi_*)[2])$.  Recall that $\calS$
denotes the coarse moduli space of supersingular genus-$2$ curves,
and in Section~\ref{S:basic3} we gave an isomorphism $\calS\to\A^1$.

\begin{theorem}
\label{T:moduli}
Let notation be as above.
\begin{enumerate}
\item The map $\kappa\colon \calA\to\calB$ is an isomorphism.
\item The function that sends a cover $\varphi_{b,c}$ to
      $c^2/b^3$ induces an isomorphism $\calA\to\A^1\setminus\{0\}$.
\item Under the isomorphisms $\calA\to \A^1\setminus\{0\}$ and
      $\calS\to\A^1$ given above, the natural map $\calA\to\calS$ that
      sends $(C,E,\varphi)$ to $C$ is given by the function
      $$t\mapsto -\frac{(1+t^4)^5}{t^{18}}.$$
\end{enumerate}
\end{theorem}

\begin{proof}
If $C$ is a genus-$2$ curve over a field of characteristic not~$2$, there
are $35$ geometric order-$4$ subgroups of $(\Jac C)[2]$.  Of these 
subgroups, there are $15$ that are isotropic with respect to the Weil 
pairing.  Therefore, the degree of the natural map $\calB\to\calS$ is~$20$.

Let $\calC$ denote $\A^1\setminus\{0\}$, and let $\iota$ be the map from
$\calC$ to $\calA$ that sends $c$ to $C_{1,c}$.  
Lemma~\ref{L:IsomorphicCovers} shows that $c$ and $-c$ have the same image
in $\calA$, so the degree of $\iota$ is at least $2$.  It is easy to
compute that 
\begin{equation}
\label{EQ:invariant}
I_{C_{1,c}} = -\frac{(1 + c^8)^5}{c^{36}},
\end{equation}
so the composition $\calC\to\calA\to\calS$ has degree $40$ and the map
$\calA\to\calS$ has degree at most $20$.  But $\calA\to\calS$ factors
through the degree-$20$ map $\calB\to\calS$, so we deduce that
$\kappa\colon\calA\to\calB$ is an isomorphism and that $\iota$ has 
degree $2$.  In particular, this proves statement~(1).

We also see that $\calA$ is isomorphic to the quotient of $\calC$ by 
$\{\pm1\}$.  This quotient is itself isomorphic to $\A^1\setminus\{0\}$,
and the isomorphism can be chosen so that the point on $\calA$ 
corresponding to $\varphi_{b,c}$ gets sent to the point $c^2/b^3$ of
$\A^1\setminus\{0\}$.  This proves statement~(2).

Statement~(3) follows from equation~\eqref{EQ:invariant}.
\end{proof}

\begin{remark}
Suppose $C$ is a genus-$2$ curve defined by $y^2 = f$, where $f$ is a 
degree-$6$ polynomial.  Suppose $g$ is a cubic factor of $f$.  Then we 
can associate to $g$ the subgroup of $(\Jac C)[2]$ whose geometric points
are the divisor classes that can be represented by divisors of the form
$P-Q$, where $P$ and $Q$ are Weierstra\ss\ points whose $x$-co\"ordinates
are roots of $g$.  The Weil pairing on two such divisors is given by the
cardinality (modulo $2$) of the intersection of their supports, so the
subgroup associated to $g$ is not isotropic.  Since there are $20$ cubic 
factors of $f$ (over the algebraic closure) and $20$ non-isotropic
geometric order-$4$ subgroups of $(\Jac C)[2]$, we can view the moduli 
space $\calB$ as the space of pairs $(f,g)$, where $y^2 = f$ is a
supersingular genus-$2$ curve and $g$ is a cubic factor of~$f$.
\end{remark}

\section{Supersingular elliptic curves over finite fields
of characteristic three}
\label{S:SSEC}

In this section we list the different isomorphism classes of supersingular 
elliptic curves over finite fields of characteristic~$3$.  We specify how
the different isomorphism classes are Galois twists of one another, and we
compute their traces of Frobenius.  This information will enable us to 
compute the Weil polynomials of supersingular genus-$2$ curves in 
Section~\ref{S:WeilPolys3}.

Let $E$ be the elliptic curve $y^2 = x^3 - x$ over $\GF_3$.  For every 
finite extension $\GF_q$ of~$\GF_3$, we will express all the isomorphism
classes of  supersingular elliptic curves over $\GF_q$ as twists of~$E$.

Let $i\in\GF_9$ be a square root of $-1$.  Let $\iota, \omega$, and $\pi$
be the following endomorphisms of $E_{\GF_9}$:
\begin{align*}
\iota &\colon (x,y)\mapsto(-x,iy)\\
\omega&\colon (x,y)\mapsto(x-1,y)\\
\pi   &\colon (x,y)\mapsto(x^3,y^3).
\end{align*}
One can check that the following relations hold in $\End E_{\GF_9}$:
$$
\iota\omega = \omega^2\iota, \qquad
\iota\pi = -\pi\iota, \qquad
\omega\pi = \pi\omega, \qquad
\pi = 1 + 2\omega.$$
The automorphism group $A$ of $E_{\GFbar_3}$ is the group of order $12$
generated by $\iota$ and $\omega$.

Let $q = 3^d$ be a power of $3$. The twists of $E_{\GF_q}$ are catalogued
by the pointed cohomology set $H^1(\Gal(\GFbar_q/\GF_q), A)$.  We identify 
a cocyle $\Gal(\GFbar_q/\GF_q)\to A$ with the image of the $q$-power 
Frobenius automorphism in $A$; using this identification, we view a 
cohomology class as a set of elements of $A$.  We let the reader verify
that when $d$ is odd (so that the Galois action on $A$ is nontrivial) we have
$$H^1(\Gal(\GFbar_q/\GF_q), A) = 
\left\{ \{\pm1\}, \{\omega,-\omega^2\},\{-\omega,\omega^2\},
\{\pm\iota,\pm\iota\omega,\pm\iota\omega^2\}\right\},$$
and when $d$ is even (so the Galois action on $A$ is trivial) we have
\begin{multline*}
H^1(\Gal(\GFbar_q/\GF_q), A) = \\
\left\{ \{1\}, \{-1\},\{\omega,\omega^2\},\{-\omega,-\omega^2\},
\{\iota,\iota\omega,\iota\omega^2\},\{-\iota,-\iota\omega,-\iota\omega^2\}
\right\}.
\end{multline*}

Every supersingular elliptic curve $E'$ over $\GF_q$ can be written in the
form $y^2 = x^3 - bx + c$ with $b\neq 0$, and every curve of this form is
supersingular.  In Table~\ref{Table:SSECodd} we enumerate the isomorphism 
classes of the twists of $E$ in the case that $d$ is odd, and we show how
one can determine the cohomology class corresponding to a supersingular
elliptic curve $E'$ from the coefficients $b$ and $c$ of its defining 
equation.  The verification of the table is a straightforward exercise; 
we simply note that the Frobenius endomorphism of the twist of $E$ by the
cohomology class containing an automorphism $\zeta$ can be taken to be 
$\zeta\pi^d$, where $\pi$ is the $\GF_3$-Frobenius of $E$ (as above).
The table refers to the quadratic character 
$\chi_{q,2}\colon \GF_q^*\to\{\pm 1\}$ defined by 
$\chi_{q,2}(x) = x^{(q-1)/2}$ and to the absolute trace function
$\Tr\colon \GF_q\to\GF_3$.  The table also uses the notation $s(b)$ to
denote the quantity $b^{(3-q)/4}$; since $q\equiv 3\bmod 8$, we see that
$s(b)$ is always a square, and when $b$ is a square, $s(b)$ is the unique 
square root of $b$ that is itself a square.

\begin{table}
\begin{center}
\renewcommand{\arraystretch}{1.1}
\begin{tabular}{|l|l||c|c|c|}
\hline
Condition on $b$ & Condition on $c$ & $\#\Aut E'$ & Cohomology class & Trace\\
\hline\hline
$\chi_{q,2}(b) =  1$ & $\Tr c/s(b)^3 = 0$ & $6$ &$ \{\pm 1\}$            & $0$                     \\
\hline
$\chi_{q,2}(b) =  1$ & $\Tr c/s(b)^3 = 1$ & $6$ & $\{\omega,-\omega^2\}$ & \hfill $(-3)^{(d+1)/2}$    \\
\hline
$\chi_{q,2}(b) =  1$ & $\Tr c/s(b)^3 = 2$ & $6$ & $\{-\omega,\omega^2\}$ & \hfill $-(-3)^{(d+1)/2}$   \\
\hline
$\chi_{q,2}(b) = -1$ & (no condition)     & $2$ & $\{\pm\iota,\pm\iota\omega,\pm\iota\omega^2\}$ & $0$\\
\hline
\end{tabular}
\end{center} 
\vspace{1ex}
\caption{The twists $E'\colon y^2 = x^3 - bx + c$ of the elliptic curve
$y^2 = x^3 - x$ over $\GF_{3^d}$, when $d$ is odd.
}
\label{Table:SSECodd} 
\end{table}

In Table~\ref{Table:SSECeven} we enumerate the isomorphism classes of the
twists of $E$ in the case that $d$ is even.  The table refers to the 
quartic character $\chi_{q,4}\colon \GF_q^*\to\{\pm 1,\pm i\}$ defined 
by $\chi_{q,2}(x) = x^{(q-1)/4}$ and to the absolute trace function 
$\Tr\colon \GF_q\to\GF_3$.  The conditions that involve $b^{3/2}$ do not
depend on the choice of the square root of~$b$.

\begin{table}
\begin{center}
\renewcommand{\arraystretch}{1.1}
\begin{tabular}{|l|l||c|c|c|}
\hline
Condition on $b$     & Condition on $c$       & $\#\Aut E'$ & Cohomology class                          & Trace\\
\hline\hline
$\chi_{q,4}(b) =  1$ & $\Tr c/b^{3/2}   =  0$ & $12$        & $\{1\}$                                   & \hfill $2(-3)^{d/2}$ \\
\cline{2-5}
                     & $\Tr c/b^{3/2} \neq 0$ & $6$         & $\{ \omega, \omega^2\}$                   & \hfill $-(-3)^{d/2}$\\
\hline
$\chi_{q,4}(b) = -1$ & $\Tr c/b^{3/2}   =  0$ & $12$        &$\{-1\}$                                   & \hfill $-2(-3)^{d/2}$ \\
\cline{2-5}
                     & $\Tr c/b^{3/2} \neq 0$ & $6$         & $\{-\omega,-\omega^2\}$                   & \hfill $(-3)^{d/2}$\\
\hline
$\chi_{q,4}(b) =  i$ & (no condition)         & $4$         & $\{-\iota,-\iota\omega,-\iota\omega^2\}$  & $0$\\
\hline
$\chi_{q,4}(b) = -i$ & (no condition)         & $4$         & $\{\iota,\iota\omega,\iota\omega^2\}$     & $0$\\
\hline
\end{tabular}
\end{center} 
\vspace{1ex}
\caption{The twists $E'\colon y^2 = x^3 - bx + c$ of the elliptic curve
$y^2 = x^3 - x$ over $\GF_{3^d}$, when $d$ is even.
}
\label{Table:SSECeven} 
\end{table}

\section{Weil polynomials of supersingular genus-two curves 
in characteristic three}
\label{S:WeilPolys3}

In this section we determine the polynomials that occur as Weil polynomials
of supersingular genus-$2$ curves over finite fields of characteristic~$3$.
In fact, for every polynomial that occurs as the Weil polynomial of a 
supersingular genus-$2$ curve over~$\GF_q$, we show how to quickly construct
an explicit curve with this Weil polynomial.

To begin, let us record for convenience a list of the Weil polynomials of
supersingular abelian surfaces over the finite fields of characteristic~$3$.
\begin{lemma}
\label{L:SSAS}
Let $q = 3^d$ be a power of $3$.  If $d$ is odd, then the polynomials
that occur as the Weil polynomials of supersingular abelian surfaces
over $\GF_q$ are{\rm:}
\begin{enumerate}
\item $(x^2 - sx + q)(x^2 - tx + q)$ for all $s,t\in\{0,\pm\sqrt{3q}\}${\rm;}
\item $x^4 + q^2${\rm;}
\item $x^4 + qx^2 + q^2${\rm;}
\item $x^4 - 2qx^2 + q^2$.
\end{enumerate}
If $d$ is even, the polynomials that occur as the Weil polynomials of 
supersingular abelian surfaces over $\GF_q$ are{\rm:}
\begin{enumerate}
\item $(x^2 - sx + q)(x^2 - tx + q)$ for all 
      $s,t\in\{0,\pm\sqrt{q},\pm2\sqrt{q}\}${\rm;}
\item $x^4 + q^2${\rm;}
\item $x^4 - qx^2 + q^2${\rm;}
\item $x^4 - sx^3 + qx^2 - sq x+ q^2$ for all $s\in\{\pm\sqrt{q}\}$.
\end{enumerate}
\end{lemma}

\begin{proof}
We already know the Weil polynomials of the supersingular elliptic curves
over $\GF_q$, and on each of the lists in the lemma, item (1) gives the
Weil polynomials of the products of these curves.  The remaining polynomials
on each list correspond to simple abelian surfaces; this follows from
\cite[Cor.~2.8, Thm.~2.9]{MaisnerNart:EM}. The cited results also show
that the lists are complete.
\end{proof}

To prove Theorem~\ref{T:intro}, then, we must go through the lists in
Lemma~\ref{L:SSAS} and determine which entries come from genus-$2$ curves.  

\begin{proof}[Proof of Theorem~{\rm\ref{T:intro}}]
We begin with the case when the degree of $\GF_q$ over $\GF_3$ is even.
The result when $q=9$ can be obtained by direct computation, so we will
assume that $q > 9$.  We begin by providing, for every Weil polynomial in 
the second list in Theorem~\ref{T:intro}, a curve with that Weil polynomial.

\emph{Items {\rm(1)}, {\rm(2)}, and {\rm(3):} The split polynomials}.
For each of the polynomials listed in items (1), (2), and (3) we can find
a curve $C_{b,c}$ with that Weil polynomial.  This is clear from 
Table~\ref{Table:SSspliteven}, in which we give simple conditions on $b$ 
and $c$ that determine the Weil polynomial of the curve $C_{b,c}$.
It is easy to use Lemma~\ref{L:SSsplitting} and Table~\ref{Table:SSECeven} 
to verify that the entries of Table~\ref{Table:SSspliteven} are correct.

\begin{table}
\begin{center}
\renewcommand{\arraystretch}{1.1}
\begin{tabular}{|l|l|c|c|}
\hline
Condition on $b$         & Condition on $c$                                 & $s$            & $t$             \\
\hline\hline
                         & $\Tr_{\GF_q/\GF_9} c/b^{3/2}  =  0$              & \hfill $ 2(-3)^{d/2}$ & \hfill $ 2(-3)^{d/2}$  \\
\cline{2-4}
$\chi_{q,4}(b) =  1$     & $\Tr_{\GF_q/\GF_9} c/b^{3/2}  = \pm1\pm i$       & \hfill $ -(-3)^{d/2}$ & \hfill $ -(-3)^{d/2}$  \\
\cline{2-4}
                         & $\Tr_{\GF_q/\GF_9} c/b^{3/2} \in\{\pm1,\pm i\}$  & \hfill $ 2(-3)^{d/2}$ & \hfill $ -(-3)^{d/2}$  \\
\hline
                         & $\Tr_{\GF_q/\GF_9} c/b^{3/2}  =  0$              & \hfill $-2(-3)^{d/2}$ & \hfill $-2(-3)^{d/2}$  \\
\cline{2-4}
$\chi_{q,4}(b) = -1$     & $\Tr_{\GF_q/\GF_9} c/b^{3/2}  = \pm1\pm i$       & \hfill $  (-3)^{d/2}$ & \hfill $  (-3)^{d/2}$  \\
\cline{2-4}
                         & $\Tr_{\GF_q/\GF_9} c/b^{3/2} \in\{\pm1,\pm i\}$  & \hfill $-2(-3)^{d/2}$ & \hfill $  (-3)^{d/2}$  \\
\hline
$\chi_{q,4}(b) =  \pm i$ & (no condition)                                   & $           0$ & $           0$  \\
\hline
\end{tabular}
\end{center} 
\vspace{1ex}
\caption{The conditions on $b$ and $c$ that determine the Weil polynomial 
$(x^2 - s x + q)(x^2 - t x + q)$ of the curve $C_{b,c}$, when $q = 3^d$ and
$d$ is even.  The function $\chi_{q,4}$ is the character defined in 
Section~\ref{S:covers}.
}
\label{Table:SSspliteven} 
\end{table}

\emph{Item {\rm(4):} The polynomial $x^4 + q^2$}.
Let $a\in\GF_q$ be a nonsquare, and let $C$ be the curve
$$y^2 = z^6 + z^3 + az + a^3 + 1$$
over $\GF_q$.  Over $\GF_{q^2}$ we can write $C$ as
$$y^2 = (z^3 - \sqrt{-a}z^2 + az - a\sqrt{-a} - 1)
        (z^3 + \sqrt{-a}z^2 + az + a\sqrt{-a} - 1)$$
so that $C_{\GF_{q^2}}$ is isomorphic to $C_{\sqrt{-a},1}$.  It is easy to
check that the isogeny 
$$\Jac C_{\GF_{q^2}} \to E_{\sqrt{-a},1} \times E_{-\sqrt{-a},1}$$
from Lemma~\ref{L:SSsplitting} descends to give an isogeny from $\Jac C$
to the restriction of scalars of $E_{\sqrt{-a},1}$.  We see from 
Table~\ref{Table:SSECeven} that the Weil polynomial of $E_{\sqrt{-a},1}$ 
is $x^2 + q^2$ (because the square root of a nonsquare in $\GF_q$ is a
nonsquare in $\GF_{q^2}$, when $q$ is a power of $9$, because the group
$\GF_{q^2}^* / \GF_q^*$ has order $2$ mod $4$), so the Weil polynomial 
for $C$ is $x^4 + q^2$. 

\emph{Item {\rm(5):} The polynomials $x^4 - sx^3 + qx^2 - sqx + q^2$,
where $s^2 = q$}.  
Let $f= x^4 - sx^3 + qx^2 - sqx + q^2$, where $s$ is a square root of~$q$.
The main result of~\cite{HMNR} shows that the isogeny class of abelian
surfaces corresponding to~$f$ contains a principally polarized variety~$A$.
We see from~\cite[Table~1, p.~325]{MaisnerNart:EM} that the simple variety
$A$ remains simple over $\GF_{q^2}$, and~\cite[Thm.~4.1]{MaisnerNart:EM}
then shows that $A$ is the Jacobian of a curve over~$\GF_q$.

(One can also construct explicit curves with these Weil polynomials as
follows.  Note that the map $\GF_q\to\GF_q$ given by $z\mapsto z^5 - z^2$
is not surjective, because $0$ and $1$ both get sent to $0$.  Therefore 
there is a nonzero $c\in\GF_q$ such that the polynomial $g = z^5 - z^2 - c$
has no roots in~$\GF_q$.  The discriminant of $g$ is $-c^4$, which is a
square in~$\GF_q$, so a result of Pellet~\cite{Pellet:CRASP} shows that $g$
has an even number of irreducible factors of even degree.  Since $g$ has no 
linear factors, it follows that $g$ is irreducible.  Then $\GF_{q^5}$ is
the smallest extension of $\GF_q$ over which the Jacobian of the 
supersingular curve $C\colon y^2 = z(z^5 - z^2 - c)$ has rational 
$2$-torsion points, so the Weil polynomial of $C$ is not any of the 
polynomials listed in items (1) through (3) of the second list in 
Lemma~\ref{L:SSAS}.  The only possibility remaining is that $C$ and its 
quadratic twist give us the two Weil polynomials 
$x^4 - sx^3 + qx^2 - sqx + q^2$, where $s^2 = q$.)

Of the polynomials given in the second list of Lemma~\ref{L:SSAS}, we have
shown that those that appear in the second list of Theorem~\ref{T:intro}
do occur as the Weil polynomials of Jacobians; now we must show that the
remaining polynomials from the second list of Lemma~\ref{L:SSAS} do not 
come from Jacobians.  We begin with the split isogeny classes.

In Table~\ref{Table:SSsplitbad} we list the split Weil polynomials over 
$\GF_q$ that we must show do not occur.  Suppose $C$ is a curve over
$\GF_q$ with one of these Weil polynomials.  Write $C$ as $y^2 = f$ for 
some sextic polynomial $f \in \GF_q[z]$.  If $f$ has a cubic factor, then 
Theorems~\ref{T:cubicfactors} and~\ref{T:SSnormalform}, together with 
Lemma~\ref{L:SSsplitting} and Table~\ref{Table:SSECeven}, show that the 
Weil polynomial for $C$ cannot factor in the way we have assumed.  (The
key fact to notice is that the quartic character $\chi_{q,4}$ takes the
same value on $b$ and on $-b$.)  Therefore, the polynomial $f$ is either
\begin{itemize}
\item an irreducible sextic;
\item the product of a linear polynomial with an irreducible quintic;
\item the product of a quadratic polynomial with an irreducible quartic; or
\item the product of three irreducible quadratics.
\end{itemize}
Corresponding to these factorizations of $f$, we obtain information
about the reduction of the Weil polynomial of $C$ modulo $2$.  In the
four cases listed above, the Weil polynomial modulo $2$ is
\begin{itemize}
\item $(x^2 + x + 1)^2$;
\item $(x^4 + x^3 + x^2 + x + 1)$;
\item $(x+1)^4$; and
\item $(x+1)^4$.
\end{itemize}

\begin{table}
\begin{center}
\renewcommand{\arraystretch}{1.1}
\begin{tabular}{|c|c||r|r|}
\hline
$s$                   & $t$                   & $S$       & $T$      \\
\hline\hline
\hfill $-2(-3)^{d/2}$ & \hfill $ -(-3)^{d/2}$ & $ 2 q^d$  & $ - q^d$ \\
\hline
\hfill $ -(-3)^{d/2}$ &        $           0$ & $ - q^d$  & $-2 q^d$ \\
\hline
       $           0$ & \hfill $  (-3)^{d/2}$ & $-2 q^d$  & $ - q^d$ \\
\hline
\hfill $  (-3)^{d/2}$ & \hfill $ 2(-3)^{d/2}$ & $ - q^d$  & $ 2 q^d$ \\
\hline
\hfill $-2(-3)^{d/2}$ &        $           0$ & $ 2 q^d$  & $-2 q^d$ \\
\hline
       $           0$ & \hfill $ 2(-3)^{d/2}$ & $-2 q^d$  & $ 2 q^d$ \\
\hline
\hfill $ -(-3)^{d/2}$ & \hfill $  (-3)^{d/2}$ & $ - q^d$  & $ - q^d$ \\
\hline
\hfill $-2(-3)^{d/2}$ & \hfill $ 2(-3)^{d/2}$ & $ 2 q^d$  & $ 2 q^d$ \\
\hline
\end{tabular}
\end{center} 
\vspace{1ex}
\caption{Certain split Weil polynomials 
$(x^2 - sx + q)(x^2 - tx + q)$ over~$\GF_q$, and the
corresponding Weil polynomials
$(x^2 - Sx + q)(x^2 - Tx + q)$ over~$\GF_{q^2}$.  Here
$q = 3^d$ for an even number $d$.
}
\label{Table:SSsplitbad} 
\end{table}

Note that the first four Weil polynomials in Table~\ref{Table:SSsplitbad}
do not factor modulo $2$ in any of these ways, so these first four
polynomials do not occur as the Weil polynomials of Jacobians over $\GF_q$.

For the next two entries in the table, we notice that the corresponding Weil
polynomials over $\GF_{q^2}$ are both $(x^2 - 2q x + q^2)(x^2 + 2q x + q^2)$.
We see that $f$ is therefore not the product of a linear polynomial and an
irreducible quintic.  Looking at the other possible factorizations of $f$
over $\GF_q$, we see that $f$ must have a cubic factor over $\GF_{q^2}$.
Once again, Theorem~\ref{T:cubicfactors}, Theorem~\ref{T:SSnormalform}, 
Lemma~\ref{L:SSsplitting}, and Table~\ref{Table:SSECeven} show that this 
Weil polynomial over $\GF_{q^2}$ cannot occur.

For the penultimate entry in Table~\ref{Table:SSsplitbad}, we again use 
Theorem~\ref{T:cubicfactors}, Theorem~\ref{T:SSnormalform}, 
Lemma~\ref{L:SSsplitting}, and Table~\ref{Table:SSECeven}.  We find that
over $\GF_{q^2}$, our hypothetical curve $C$ must be isomorphic to $C_{b,c}$
for some $b,c\in\GF_{q^2}$ where $b$ is a fourth power and where the 
absolute traces of $c/b^{3/2}$ and of $c/(-b)^{3/2}$ are both nonzero.
Using Lemma~\ref{L:IsomorphicCovers}, we can assume that $b= 1$, so that
$c\in\GF_{q^2}$ now has the property that $c$ and $ic$ have nonzero trace
to~$\GF_3$.  It must certainly be the case that $c$ has nonzero trace 
to~$\GF_9$.  Now, we know from equation~\eqref{EQ:invariant} and from the
fact that $C$ is defined over $\GF_q$ that we must have
$$\frac{(c^8+1)^5}{c^{36}}\in \GF_q.$$
Applying Lemma~\ref{L:HatRabbitEven} (below) with $r = 9$, we find that $c$
must actually lie in $\GF_q$.  From this it follows that $C$ is a twist of
the curve $C_{1,c}$ over $\GF_q$.  But again using Lemma~\ref{L:SSsplitting}
and Table~\ref{Table:SSECeven}, we see that $C_{1,c}$ has Weil polynomial
$(x^2 + 3^{d/2} x + q)^2$, so the Weil polynomial of $C$ must be either
$(x^2 + 3^{d/2} x + q)^2$ or $(x^2 - 3^{d/2} x + q)^2$.
This contradicts our assumption that $C$ has Weil polynomial
$(x^2 + 3^{d/2} x + q)(x^2 - 3^{d/2} x + q)$.

We turn to the final entry in Table~\ref{Table:SSsplitbad}.  Write our 
hypothetical curve $C$ as $d y^2 = f$ for some nonzero $d$ in $\GF_q$ and 
some degree-$6$ polynomial~$f$.  From what we have already noted, we know 
that $f$ is either the product of three irreducible quadratics, or is the 
product of an irreducible quartic with a quadratic.  In either case $f$ has
a cubic factor over $\GF_{q^2}$. Theorem~\ref{T:cubicfactors},
Theorem~\ref{T:SSnormalform}, Lemma~\ref{L:SSsplitting}, and 
Table~\ref{Table:SSECeven} show that $C_{\GF_{q^2}}$ is isomorphic to 
$C_{b,c}$ for some $b,c\in\GF_{q^2}^*$, where $b$ is a fourth power and
where the trace of $c/b^{3/2}$ from $\GF_{q^2}$ to $\GF_9$ is equal to~$0$.  
Without loss of generality, we may assume that $b=1$. 
Equation~\eqref{EQ:invariant} shows that $(1+c^8)^5/c^{36}$ lies 
in~$\GF_{q^2}$ (indeed, it lies in~$\GF_q$), and applying 
Lemma~\ref{L:HatRabbitEvenBis} (below) with $r=9$ and $Q=q^2$ shows that 
all of the geometric points of $\calA$ lying above the point $P_C$ of 
$\calS$ corresponding to $C$ are rational over $\GF_{q^2}$.  Therefore all
of the geometric points of $\calB$ lying above $P_C$ are rational 
over~$\GF_{q^2}$, and it follows that all of the geometric cubic factors 
of $f$ are rational over $\GF_{q^2}$; in other words, $f$ splits completely
over $\GF_{q^2}$.  Hence we see that $f$ must be the product of three 
irreducible quadratics over $\GF_q$.

At the beginning of the proof of Theorem~\ref{T:SSfacts1} we showed that
every supersingular genus-$2$ curve over $\GF_q$ can be written in the
form $dy^2 = z^6 + c_3 z^3 + c_1 z + c_0$, so we may assume that our
polynomial $f$ has the form of the right-hand side of this equation.
Then the discriminant of $f$ is $-c_1^6$, and since $\GF_q$ is an 
even-degree extension of $\GF_3$, the discriminant is a square.  But a 
polynomial over a finite field with a square discriminant has an even
number of irreducible even-degree factors~\cite{Pellet:CRASP}, so $f$ 
cannot be the product of three irreducible quadratics.  This contradiction 
shows that there is no supersingular curve over $\GF_q$ whose Weil polynomial
is given by the final entry in Table~\ref{Table:SSsplitbad}.

Finally, we must show that there is no supersingular genus-$2$ curve with
Weil polynomial $x^4 - qx^2 + q^2$.  If there were such a curve $C$ over 
$\GF_q$, then reasoning as above shows that $C$ must be given by an 
equation $y^2 = \text{(sextic)}$, which shows that $C_{\GF_{q^2}}$ has a
degree-$3$ map to an elliptic curve.  In fact, since the quadratic twist
$D$ of $C_{\GF_{q^2}}$ has Weil polynomial $(x^2 + qx + q^2)^2$, we find
from Table~\ref{Table:SSspliteven} that $D$ is isomorphic to $C_{1,c}$ for
some $c\in\GF_{q^2}$ whose trace to $\GF_9$ is $\pm1\pm i$.  But since $D$
is a twist of a curve over $\GF_q$, we know that 
$$I_D = -\frac{(c^8+1)^5}{c^{36}}$$
is an element of $\GF_q$.  Then Lemma~\ref{L:HatRabbitEven} shows that $c$
must be an element of $\GF_q$, so our original curve $C$ is a twist of a 
curve $C_{1,c}$.  But $C_{1,c}$ has a split Jacobian and $C$ does not, and
this contradiction shows that there is no curve over $\GF_q$ with Weil 
polynomial $x^4 - qx^2 + q^2$.

This completes the proof of Theorem~\ref{T:intro} in the case that $d$ is
even.  We now turn to the case in which $d$ is odd.  The results when $q=3$
can be obtained by enumerating all of the supersingular genus-$2$ curves 
over $\GF_3$ and computing their Weil polynomials, so we will assume 
that $q > 3$.

Let us begin by giving examples of curves whose Weil polynomials comprise
all the polynomials given in the first list of Theorem~\ref{T:intro}.

\emph{Items {\rm(1)} and {\rm(2):} The split polynomials.}
Let $c_0$, $c_1$, and $c_2$ be nonzero elements of $\GF_q$ whose traces to
$\GF_3$ are $0$, $1$, and $2$, respectively.  Using Lemma~\ref{L:SSsplitting}
and Table~\ref{Table:SSECodd}, it is easy to check that the curve $C_{1,c_1}$
has Weil polynomial
$$(x^2 + q)(x^2 - (-1)^{(d+1)/2} \sqrt{3q} x + q),$$
that the curve $C_{1,c_2}$ has Weil polynomial
$$(x^2 + q)(x^2 + (-1)^{(d+1)/2} \sqrt{3q} x + q),$$
and that the curve $C_{1,c_0}$ has Weil polynomial
$$(x^2 + q)^2.$$

\emph{Item {\rm(3):} The polynomial $x^4 + q^2$}.
We calculate that the supersingular curve  $y^2=z^5+1$ over $\GF_3$ has 
Weil polynomial $x^4 + 9$.  It follows that the Weil polynomial for this
curve over $\GF_q$ is $x^4 + q^2$ whenever $q$ is an odd power of~$3$.

\emph{Item {\rm(4):} The polynomial $x^4 + qx^2 + q^2$}.
Let $a$ be an element of $\GF_q$ whose absolute trace is nonzero.
Let $i\in\GF_9$ be a square root of $-1$, let $b$ be the element $i a^2$
of $\GF_{q^2}$, and let $c = -a^4$.  Let $C$ be the curve 
$$y^2 = -a^4(z^6 + a^8 z^3 + a^{12} z + a^{16} + a^{12})$$
over $\GF_q$.  Note that over $\GF_{q^2}$ we have
$$
z^6 + a^8 z^3 + a^{12} z + a^{16} + a^{12} 
    = (z^3 -ia^2 z^2 + a^4 z -ia^6-a^8)
           (z^3 +ia^2 z^2 + a^4 z +ia^6-a^8)$$
so that $C_{\GF_{q^2}} \cong C_{b,c}$.  If we identify $C_{b,c}$ with 
$C_{\GF_{q^2}}$ in the obvious way then the Galois group of $\GF_{q^2}$ over
$\GF_q$ interchanges the two morphisms $\varphi_{b,c}\colon C\to E_{b,c}$ and
$\varphi_{b,c}'\colon C\to E_{-b,c}$.  It follows that the isogeny 
$$\Jac C_{\GF_{q^2}} \to E_{b,c} \times E_{-b,c}$$
over $\GF_{q^2}$ descends to give an isogeny 
$$\Jac C \to \Res_{\GF_{q^2}/\GF_q} E_{b,c}$$
from the Jacobian of $C$ to the restriction of scalars of $E_{b,c}$.  Since
$i$ is a square but not a fourth power in $\GF_9$, it is easy to see that 
$\chi_{q^2,4}(i a^2) = -1$.  It is also easy to check that
\begin{align*}
\Tr_{\GF_{q^2}/\GF_3} c/b^{3/2}
 &= \Tr_{\GF_{q^2}/\GF_3} (1-i) a
 = \Tr_{\GF_9/\GF_3} \Tr_{\GF_{q^2}/\GF_9} (1-i)a\\
 &= \Tr_{\GF_9/\GF_3} (1-i)\Tr_{\GF_{q^2}/\GF_9} a
 = \Tr_{\GF_9/\GF_3} (1-i)\Tr_{\GF_{q}/\GF_3} a\\
 &= - \Tr_{\GF_{q}/\GF_3} a \neq 0,
\end{align*}
so by Table~\ref{Table:SSECeven} we see that the Weil polynomial for 
$E_{b,c}$ is $x^2 + qx + q^2$.  It follows that the Weil polynomial for 
$\Res_{\GF_{q^2}/\GF_q} E_{b,c}$ is $x^4 + qx^2 + q^2$, so $\Jac C$ has
this Weil polynomial as well.

(Note that the quadratic twist of this $C$ also has Weil polynomials
$x^4 + qx^2 + q^2$.  By rescaling the variables in the equation for the 
twist, we find that 
$$y^2 =z^6 + a^2 z^3 + a^2 z + a^4 + 1$$
has Weil polynomial $x^4 + qx^2 + q^2$ whenever $a$ is an element of 
$\GF_q$ whose absolute trace is nonzero.)

\emph{Item {\rm(5):} The polynomial $x^4 - 2qx^2 + q^2$}.
We use a similar construction to produce a curve with Weil polynomial
$x^4 - 2qx^2 + q^2$.  Let $a\in\GF_q$ be a nonzero element whose absolute
trace is $0$, let $i\in \GF_9$ be as above, and again take 
$b = ia^2\in\GF_{q^2}$ and $c = -a^4$.  We again find that the curve $C$ 
over $\GF_q$ defined by 
$$y^2 = z^6 + a^2 z^3 + a^2 z + a^4 + 1$$
has Jacobian isogenous to the restriction of scalars of $E_{b,c}$, but now
Table~\ref{Table:SSECeven} shows that the Weil polynomial of $E_{b,c}$ is
$x^2 - 2qx + q^2$.  It follows that the Weil polynomial of $C$ is
$x^4 - 2qx^2 + q^2$.

To complete the proof of Theorem~\ref{T:intro}, we must show that no curve 
over $\GF_q$ has Weil polynomial
$$(x^2 -\eps_1\sqrt{3q} + q) (x^2 -\eps_2\sqrt{3q} + q)$$
for any choice of $\eps_1, \eps_2 \in \{1,-1\}$.  Suppose, to obtain a 
contradiction, that $C$ is a curve over $\GF_q$ with such a Weil polynomial.
Write $C$ as $y^2 = f$ for a degree-$6$ polynomial $f\in\GF_q[z]$.   Note
that $f$ is not the product of a linear polynomial with an irreducible
quintic, because if this were the case the minimal polynomial of Frobenius
on the $2$-torsion of $C$ would be $x^4 + x^3 + x^2 + x + 1$.  By
enumerating the other possible factorization for $f$, we see that over 
$\GF_{q^2}$, the polynomial $f$ can be written a product of two 
(not-necessarily-irreducible) cubics.  Therefore, by 
Theorems~\ref{T:cubicfactors} and~\ref{T:SSnormalform}, $C_{\GF_{q^2}}$ 
is isomorphic to $C_{b,c}$ for some $b,c\in\GF_{q^2}^*$.

A simple computation shows that the Weil polynomial for $C_{\GF_{q^2}}$ is
$$(x^2 - qx + q^2)^2,$$
so Lemma~\ref{L:SSsplitting} and Table~\ref{Table:SSECeven} show that 
$\chi_{q^2,4}(b) = \chi_{q^2,4}(-b) = 1$ and that the traces of $c/b^{3/2}$
and of $c/(-b)^{3/2}$ from $\GF_{q^2}$ to $\GF_3$ are nonzero.  Using 
Lemma~\ref{L:IsomorphicCovers}, we can assume without loss of generality
that $b = 1$.  Then we know that $\Tr_{\GF_{q^2}/\GF_3} c \neq 0$, and in
particular $\Tr_{\GF_{q^2}/\GF_9} c$ is nonzero.

Equation~\eqref{EQ:invariant} tells us that 
$$I_C = -\frac{(c^8+1)^5}{c^{36}}.$$
In particular, we note that $I_C\neq 0$, because the roots of $z^8 + 1$
generate $\GF_{81}$, and $\GF_{q^2}$ is not an extension of $\GF_{81}$. 
It follows from Theorem~\ref{T:SSfacts2}(c) that the geometric automorphism
group of $C$ has order $2$.  We also know that $I_C$ is an element of 
$\GF_q$, because $C$ is a curve over~$\GF_q$.

Using Lemma~\ref{L:HatRabbitOdd} (below), we see that there is an element
$\cbar\in\GF_q$ such that 
$$-\frac{(\cbar^8+1)^5}{\cbar^{36}} = \pm I_C.$$
If the equality holds with the plus sign, take $\bbar = 1$;  if it holds
with the minus sign, take $\bbar = i$, where $i$ is a square root of $-1$ 
in $\GF_9$.  Then we have $I_{C_{\bbar,\cbar}} = I_C$.

We claim that there is a curve $D$ over $\GF_q$ such that 
$D_{\GF_{q^2}}\cong C_{\bbar,\cbar}$.  If $\bbar=1$ then this is obvious.
If $\bbar = i$, then we note that $C_{\bbar,\cbar}$ is defined by the 
equation
\begin{align*}
y^2 &= \cbar(z^3 - i z^2 + z - i - \cbar^{2})(z^3 + i z^2 + z + i - \cbar^{2})\\
    &= \cbar(z^6 + \cbar^{2} z^3 + \cbar^{2} z + \cbar^{4} + 1),
\end{align*}
so we can take $D$ to be the curve over $\GF_q$ defined by this equation.

The fact that $I_D = I_C$ shows that $D$ is a twist of $C$;  since the
geometric automorphism group of $C$ has order $2$ the only nontrivial twist
of $C$ is the standard quadratic twist.  This shows that the Weil polynomial
of the curve $D$ is of the form
\begin{equation}
\label{EQ:smallfield}
(x^2 -\epsbar_1\sqrt{3q} + q) (x^2 -\epsbar_2\sqrt{3q} + q)
\end{equation}
for some $\epsbar_1,\epsbar_2\in\{1,-1\}$, and the Weil polynomial of
$D_{\GF_{q^2}}$ is 
\begin{equation}
\label{EQ:bigfield}
(x^2 - qx + q^2)^2.
\end{equation}

If $\bbar = 1$, then we see from Lemma~\ref{L:SSsplitting} and 
Table~\ref{Table:SSECodd} that the Weil polynomial of $D$ cannot be of the
form given by~\eqref{EQ:smallfield}.  If $\bbar = i$, then we note that
$\chi_{q^2,4}(i) = -1$, so that Lemma~\ref{L:SSsplitting} and
Table~\ref{Table:SSECeven} show that the Weil polynomial
of $D_{\GF_{q^2}}$ cannot be of the form given by~\eqref{EQ:bigfield}.
These contradictions shows that there can be no curve $C$ over $\GF_q$ 
with Weil polynomial 
$$(x^2 -\eps_1\sqrt{3q} + q) (x^2 -\eps_2\sqrt{3q} + q)$$
for any choice of $\eps_1, \eps_2 \in \{1,-1\}$.  
\end{proof}
        
We end by proving the three lemmas that we needed in the preceding proof.

\begin{lemma}
\label{L:HatRabbitEven}
Let $r$ be a power of an odd prime $p$, and let $q$ be a power of $r$.
If $c\in\GF_{q^2}\setminus\GF_q$ satisfies 
$(c^{r-1} + 1)^{(r+1)/2}/c^{r(r-1)/2}\in \GF_q$, then 
$\Tr_{\GF_{q^2}/\GF_r} c = 0$.
\end{lemma}

\begin{proof}
Let $G = \PSL_2(\GF_r)$, and consider the rational function
$$F = \frac{ \left( (z^r-z)^{r-1} + 1\right)^{(r+1)/2}}
           {(z^r-z)^{r(r-1)/2}} 
    = \frac{ (z^{r^2} - z)^{(r+1)/2}}
           {(z^r-z)^{(r^2+1)/2}}.$$
It is easy to verify that $F$ is fixed by the action of $G$ on $\GF_r(z)$
by fractional linear transformations.  Since the degree of $F$ is equal to
$\#G$, we know that for every $z_0\in\GFbar_r$, the roots of 
$F(z) - F(z_0)$ in $\GFbar_r$ are precisely the images of $z_0$ under the
action of $G$.  Furthermore, if $F(z_0)$ is nonzero then one can verify that 
$$\left( (z^r-z)^{r-1} + 1\right)^{(r+1)/2} - F(z_0)(z^r-z)^{r(r-1)/2}$$
is a separable polynomial, so $G$ acts faithfully on the roots of 
$F(z) -F(z_0)$.

Let $c$ be as in the statement of the lemma.  Suppose that the trace of $c$
to $\GF_r$ were not $0$.  Then the polynomial $z^r - z - c$ would have no
roots in $\GF_{q^2}$, so all of the roots of $z^r - z - c$ live in an 
extension of $\GF_{q^2}$ of degree divisible by $p$.  Let $d$ be such a 
root, and let $n = [\GF_q(d): \GF_q]$, so that $n$ is a multiple of~$2p$.
Then the hypothesis of the first statement of the lemma is that 
$F(d)\in\GF_q$.

Since $d^q$ is another root of $F(z) - F(d)$, we know there is an element
$\sigma$ of $\PSL_2(\GF_r)$ so that $d^q = \sigma(d)$.  The element 
$\sigma$ of $G$ is fixed by the action of the absolute Galois group of
$\GF_q$, so we find that $d^{q^i} = \sigma^i(d)$ for all integers $i\ge 0$.
Thus $\sigma^n$ is the smallest power of $\sigma$ that fixes $d$.  In
particular, the order of $\sigma$ is a multiple of $2p$.  But there are 
no elements of order $2p$ in $G$, because there are no involutions in $G$
that commute with a nontrivial translation.  This contradiction proves the
lemma.
\end{proof}

\begin{lemma}
\label{L:HatRabbitEvenBis}
Let $r$ be a power of an odd prime $p$, and let $Q$ be a power of $r$.
Suppose $c\in\GF_Q$ satisfies $\Tr_{\GF_Q/\GF_r} c = 0$, and let 
$e = (c^{r-1} + 1)^{(r+1)/2}/c^{r(r-1)/2}$.  Then the polynomial
$$(z^{r-1}+1)^{(r+1)/2} - e z^{r(r-1)/2}$$ 
splits completely over~$\GF_Q$.
\end{lemma}

\begin{proof}
Since $\Tr_{\GF_Q/\GF_r} = 0$, there is a $d\in\GF_Q$ with $c = d^r - d$.
Let $F$ be as in the preceding proof.  Then all of the zeroes of the 
function $F(z) - F(d) = F(z) - e$ lie in $\GF_Q$, and it follows that the
polynomial in the statement of the lemma splits completely.
\end{proof}

\begin{lemma}
\label{L:HatRabbitOdd}
Let $q$ be an odd power of $3$.  Suppose $c\in\GF_{q^2}\setminus\GF_q$
satisfies $(c^8+1)^5/c^{36}\in\GF_q$ and that $\Tr_{\GF_{q^2}/\GF_9}(c)$
is nonzero.  Then there is a $\cbar\in\GF_q$ such that
$$\frac{(\cbar^8+1)^5}{\cbar^{36}}  = \pm \frac{(c^8+1)^5}{c^{36}}.$$
\end{lemma}

\begin{proof}
We use many of the ideas from the proof of Lemma~\ref{L:HatRabbitEven}.
Let $G = \PSL_2(\GF_9)$ and let $F$ be the rational function
$$F = \frac{ \left( (z^9-z)^{8} + 1\right)^{5}}
           {(z^9-z)^{36}}
    = \frac{ (z^{81}-z)^5}
           {(z^9 - z)^{41}}.$$
Since the trace of $c$ to $\GF_9$ is nonzero, we know that the polynomial 
$z^9 - z - c$ over $\GF_{q^2}$ factors as a product of cubics.  Let $d$ be
an element of $\GFbar_q$ with $d^9 - d = c$.  Then $\GF_q(d)$ is a
degree-$6$ extension of $\GF_q$.

Let $e = F(d)$ so that 
$$e = \frac{(c^8+1)^5}{c^{36}}\in\GF_q.$$
Then $d$ is a root of the rational function $F(z) - e$ in $\GF_q(z)$, so
$d^q$ is another root of this function.  By the reasoning from the proof 
of Lemma~\ref{L:HatRabbitEven}, we see that there is an element 
$\sigma\in G$ with $d^q = \sigma(d)$.  Let $\sigmabar$ denote the image
of $\sigma$ under the action of $\Gal(\GF_9/\GF_3)$ on $G$.  Then 
$(\sigmabar\sigma)^3$ acts trivially on $d$, but no smaller power of
$(\sigmabar\sigma)$ acts trivially.   It follows that $\sigmabar\sigma$ 
has order $3$.  Enumerating the elements $\sigma$ of $G$ with this property,
we find that there are elements $\tau$ and $\rho$ of $G$ such that 
$\sigma = \taubar\rho\tau^{-1}$ and 
$$\rho \in \left\{ \left(\begin{matrix} 1&1+i\\0&1\end{matrix}\right),
\left(\begin{matrix} -1+i&1+i\\0&1+i\end{matrix}\right)\right\}.$$
Let $d' = \tau^{-1}(d)$, so that $d'$ is yet another root of $F(z) - e$, 
but now $(d')^q = \rho(d')$. 

If $\rho(d') = d' + 1 + i$ then take $\cbar = (d')^9 - d'$.  It is easy to
check that $\cbar\in\GF_q$ and that 
$$\frac{(\cbar^8+1)^5}{\cbar^{36}}  =  \frac{(c^8+1)^5}{c^{36}}.$$
If $\rho(d') = id' + 1$ then take $\cbar = (1+i)d'$. Then again we have
$\cbar\in\GF_q$, but now
$$\frac{(\cbar^8+1)^5}{\cbar^{36}}  =  -\frac{(c^8+1)^5}{c^{36}}.$$
\end{proof}



\begin{thebibliography}{99}


\bibitem{AdlemanHuang:LNM} 
{\sc Leonard M. Adleman and Ming-Deh A. Huang}:
\emph{Primality testing and abelian varieties over finite fields}, 
Lecture Notes in Math.~{\bf 1512}, Springer-Verlag, Berlin, 1992.

\bibitem{Magma}
{\sc Wieb Bosma, John Cannon, and Catherine Playoust}:
The Magma algebra system. I. The user language,
\emph{J. Symbolic Comput.}~{\bf 24} (1997) 235--265.

\bibitem{Howe:TAMS} 
{\sc Everett W. Howe}:
Principally polarized ordinary abelian varieties over finite fields,
\emph{Trans. Amer. Math. Soc.}~{\bf 347} (1995) 2361--2401.

\bibitem{Howe:JAG}
{\sc Everett W. Howe}:
Kernel of polarizations of abelian varieties over finite fields,
\emph{J. Algebraic Geom.}~{\bf 5} (1996) 583--608.

\bibitem{Howe:Texel}
{\sc Everett W. Howe}:
Isogeny classes of abelian varieties with no principal polarizations,
pp.~203--216 in: \emph{Moduli of abelian varieties}
(C. Faber, G. van der Geer and F. Oort, eds.), 
Progr. Math. {\bf 195}, Birkh\"auser, Basel, 2001.

\bibitem{Howe:CM}
{\sc Everett W. Howe}:
On the non-existence of certain curves of genus two,
\emph{Compos.Math.}~{\bf 140} (2004) 581--592.

\bibitem{HoweLauter:AIF}
{\sc E.W. Howe and K.E. Lauter}:
Improved upper bounds for the number of points on curves
over finite fields,
\emph{Ann. Inst. Fourier (Grenoble)}~{\bf 53}
(2003) 1677--1737.  

\bibitem{HoweLauter:Corr}
{\sc E.W. Howe and K.E. Lauter}:
Corrigendum:
``Improved upper bounds for the number of points on curves
over finite fields,'' in preparation, 2006.

\bibitem{HMNR} 
{\sc Everett W. Howe, Daniel Maisner, Enric Nart, and Christophe Ritzenthaler}:
Principally polarized isogeny classes of abelian surfaces
over finite fields, {\tt arXiv:math.NT/0602650}.

\bibitem{HNR}
{\sc Everett W. Howe, Enric Nart, and Christophe Ritzenthaler}:
Jacobians in isogeny classes of abelian surfaces over finite fields,
preprint, 2006.

\bibitem{IKO:CM}
{\sc T. Ibukiyama, T. Katsura, and F. Oort}:
Supersingular curves of genus two and class numbers,
\emph{Compositio Math.}~{\bf 57} (1986) 127--152.

\bibitem{Igusa:AM}
{\sc J. Igusa}:
Arithmetic variety of moduli for genus two,
\emph{Ann. of Math. (2)}~{\bf 72} (1960) 612--649. 

\bibitem{Maisner:Thesis}
{\sc Daniel Maisner},
Superficies abelianas como jacobianas de curvas en cuerpos finitos,
thesis, Universitat Aut\`onoma de Barcelona, 2004.

\bibitem{MaisnerNart:EM}
{\sc Daniel Maisner and Enric Nart with an appendix by Everett W. Howe}:
Abelian surfaces over finite fields as Jacobians,
\emph{Experiment. Math.}~{\bf 11} (2002) 321--337.

\bibitem{MaisnerNart:CJM}
{\sc D. Maisner and E. Nart}:
Zeta functions of supersingular curves of genus two,
\emph{Canad. J. Math.}, to appear.

\bibitem{McGuireVoloch:PAMS}
{\sc Gary McGuire and Jos\'e Felipe Voloch}:
Weights in codes and genus $2$ curves,
\emph{Proc. Amer. Math. Soc.}~{\bf 133} (2005) 2429--2437.

\bibitem{Pellet:CRASP}
{\sc A.-E. Pellet}:
Sur la d\'ecomposition d'une fonction enti\`ere en facteurs
irr\'eductibles suivant un module premier $p$,
\emph{C. R. Acad. Sci. Paris}~{\bf 86} (1878) 1071--1072.

\bibitem{Ruck:CM} 
{\sc Hans-Georg R\"uck}:
Abelian surfaces and Jacobian varieties over finite fields, 
\emph{Compositio Math.}~{\bf 76} (1990) 351--366.

\bibitem{Serre:notes}
{\sc Jean-Pierre Serre}:
{\it Rational points on curves over finite fields},
unpublished notes by Fernando Q. Gouv\'ea of lectures at Harvard University, 1985.

\bibitem{Yui:JA}
{\sc Noriko Yui}:
On the Jacobian varieties of hyperelliptic curves over fields
of characteristic $p > 2$,
\emph{J. Algebra}~{\bf 52} (1978) 378--410.

\bibitem{Zarhin:BSMF}
{\sc Yuri G. Zarhin}:
Non-supersingular hyperelliptic Jacobians,
\emph{Bull. Soc. Math. France}~{\bf 132} (2004) 617--634.

\bibitem{Zarhin:AGCT9}
{\sc Yuri G. Zarhin}:
Homomorphisms of abelian varieties, pp. 189--215 in:
\emph{Arithmetic, Geometry and Coding Theory (AGCT 2003)}
(Y. Aubry and G. Lachaud, eds.), 
S\'emin. Congr.~{\bf 11}, Soc. Math. France, Paris 2005.

\end{thebibliography}
\end{document}